\documentclass{amsart}
\usepackage{amsfonts}
\usepackage{amsmath,amssymb}
\usepackage{amsthm}
\usepackage{amscd}
\usepackage{graphics}
\usepackage{graphicx}

\theoremstyle{remark}{

\newtheorem{Rem}{{\rm Remark}}
\newtheorem{Prob}{{\rm Problem}}

}
\theoremstyle{plain}
{

\newtheorem{Thm}{Theorem}

}

\begin{document}
\title[On the Reeb spaces of explicit analytic functions]{A note on Reeb spaces of explicit real analytic functions}
\author{Naoki kitazawa}
\keywords{Real analytic functions and maps. Smooth maps. Graphs. Infinite graphs. Peano Continua. Reeb spaces. Reeb graphs. \\
\indent {\it \textup{2020} Mathematics Subject Classification}: Primary~26E05, 57R45, 58C05. Secondary~54C30, 54F15}

\address{Osaka Central Advanced Mathematical Institute (OCAMI) \\
3-3-138 Sugimoto, Sumiyoshi-ku Osaka 558-8585
TEL: +81-6-6605-3103
}
\email{naokikitazawa.formath@gmail.com}
\urladdr{https://naokikitazawa.github.io/NaokiKitazawa.html}
\maketitle
\begin{abstract}

{\it Reeb spaces} of smooth functions are fundamental and strong tools in understanding manifolds via smooth functions with mild critical points. They are defined as the natural spaces of all connected components of level sets. They are also important objects in related studies. Realization of graphs as Reeb spaces of smooth functions of certain nice classes is of such studies. 

In this note, we present Reeb spaces of explicit real analytic functions which are not finite graphs. 
Studies on the realization were started by Sharko, in 2006, who has studied smooth functions with critical points represented by certain elementary polynomials, and followed by a study of Masumoto and Saeki, which is on smooth functions on closed surfaces under an extended situation, and a study of Michalak, which is on Morse functions on closed manifolds. The author has contributed to this by respecting topologies of level sets, and real algebraic construction.

\end{abstract}
%【REVISE】 combinatoric ～ is → combinatorial object. It is .
%【REVISE】  such that a point is a vertex if and only if the corresponding connected component of the level set contains some singular points → whose vertex set is the set of all points containing some singular points in the corresponding connected component of the level set .
%【REVISE】 We delete "extending the result before".
\section{Introduction.}
\label{sec:1}

The {\it Reeb space} $R_c$ of a (continuous) map $c:X \rightarrow Y$ between topological spaces is the space of all connected components of all preimages $c^{-1}(q)$. Such objects are strong tools in understanding the spaces $X$ via the maps $c$ and in \cite{reeb}, they have been recognized as fundamental and strong tools in understanding the manifolds $X$ via nice real-valued functions $c:X \rightarrow \mathbb{R}$ such as Morse functions.

More precisely, the equivalence relation ${\sim}_c$ on $X$ can be defined in the following way: the relation $p_1 {\sim}_c p_2$ holds if and only if $p_1$ and $p_2$ are in a same connected component of $c^{-1}(q)$. We define $R_c:=X/{\sim}_c$. 
We also have the quotient map $q_c:X \rightarrow R_c$ with the unique continuous map $\bar{c}:R_c \rightarrow Y$ satisfying the relation $c=\bar{c} \circ q_c$. Hereafter, we consider the case of a real-valued function $c:X \rightarrow \mathbb{R}$ and in this case, we call the preimage $c^{-1}(q)$ a {\it level set} of $c$ and each connected component of each level set $c^{-1}(q)$ a {\it contour} of $c$. In the case $X$ is a differentiable manifold, a contour of $c$ containing critical points (no critical point) of $c$ is said to be {\it critical} (resp. {\it regular}).

For a smooth function $c:X \rightarrow \mathbb{R}$, in some nice situation, $R_c$ is homeomorphic to a graph. It is already shown in the case of a Morse function on a compact manifold in \cite{izar}. Later, in \cite{martinezalfaromezasarmientooliveira}, this is extended to the {\it Morse-Bott} function case. For general situations, see \cite{gelbukh1, saeki2, saeki3}, especially, \cite[Theorem 7.5]{gelbukh1}, \cite[Theorem 3.1]{saeki2}, and \cite[Theorems 2.1 and 2.8 and "2 Reeb space and its graph structure"]{saeki3}, for example. More precisely, in the case the smooth manifold $X$ has no boundary, $R_c$ is regarded as a graph whose vertex set consists of all points representing critical contours of $c$ and this is the {\it Reeb graph} of $c$. 
\begin{Prob}
\label{prob:1}
Can we construct a nice smooth function whose Reeb graph is isomorphic to a given graph?
\end{Prob}
Here, a {\it finite} graph means a $1$-dimensional finite simplicial complex. This is considered to be a $1$-dimensional finite CW complex the closure of whose $1$-cell is always homeomorphic to $D^1:=\{t \mid -1 \leq t \leq 1\}$ (we avoid graphs with loops, where we do not avoid so-called {\it multigraphs}: this rule is also applied in {\it infinite} graphs, presented later). A {\it vertex} of a graph is a $0$-cell of it and an {\it edge} of it is a $1$-cell of it.

Problem \ref{prob:1} was started by Sharko in 2006. In \cite{sharko}, Sharko has considered finite graphs of a certain class and constructed smooth functions on closed surfaces whose critical points are represented by certain elementary polynomials. This is extended to arbitrary finite graphs in \cite{masumotosaeki}. In \cite{michalak1}, Michalak has considered finite graphs of a certain class and reconstructed Morse functions whose regular contours are spheres on suitable closed manifolds of dimension at least $2$ and he has also classified Reeb graphs of given closed surfaces completely. For related studies, see also \cite{marzantowiczmichalak, michalak2} and for Morse-Bott functions, see \cite{gelbukh2, gelbukh3, martinezalfaromezasarmientooliveira} for example. The author has also contributed to related studies on reconstructing nice smooth functions by respecting not only Reeb graphs, but also topologies of (regular) contours: see \cite{kitazawa1, kitazawa2, kitazawa4} for example. In these studies, we first construct local functions corresponding to vertices and product bundles corresponding to edges and glue these local functions together to have a desired global function. We also note that in \cite{kitazawa2} the author has studied construction of {\it non-proper} functions or functions on non-compact manifolds with prescribed level sets.
Problem \ref{prob:2} is essentially started by the author in \cite{kitazawa3}.
\begin{Prob}
\label{prob:2}
Can we construct these smooth functions in the real algebraic category?
\end{Prob}
We omit exposition on the real algebraic category and related terminologies, objects and tools. For this, see \cite{bochnakcosteroy, kollar, nash}. The author has found a pioneering answer in \cite{kitazawa3}, followed by the author himself in \cite{kitazawa6}. The author has also found another answer in \cite{kitazawa5}, for example. Its main ingredient consists of construction of a map of a certain class generalizing the canonical projections of the unit spheres or the class of so-called {\it special generic} maps in real algebraic situations onto a region in the $n$-dimensional Euclidean space ${\mathbb{R}}^n$ (the $n$-dimensional {\it real affine space}) surrounded by real algebraic manifolds of dimensional $n-1$ and collapsing to finite graphs naturally respecting the projection ${\pi}_{n,1}:{\mathbb{R}}^n \rightarrow \mathbb{R}$ defined by ${\pi}_{n,1}(x_1,\cdots x_n):=x_1$ (${\pi}_{k_1,k_2}(x_1,\cdots x_{k_2},\cdots x_{k_1}):=(x_1, \cdots x_{k_2})$ with $1 \leq k_2<k_1$). By composing the projection ${\pi}_{n,1}$, we have a desired function. Here the $m$-dimensional unit sphere is $S^m:=\{(x_1,\cdots, x_{m+1}) \in {\mathbb{R}}^{m+1} \mid {\Sigma}_{j=1}^{m+1} {x_j}^2=1\}$ and by the canonical projection, the restriction of ${\pi}_{m+1,n}$, it is mapped onto the $n$-dimensional unit disk $D^n:=\{(x_1,\cdots x_n) \mid  {\Sigma}_{j=1}^{n} {x_j}^2 \leq 1\}$ ($n \leq m$). Hereafter, we use $0 \in {\mathbb{R}}^k$ for the origin.
For {\it special generic} maps, see \cite{saeki1} for example. For graphs the regions collapses to, see \cite{bodinpopescupampusorea}. Different from the differentiable category, it is difficult to construct global functions and maps, due to the non-existence of objects or facts such as a partition of the unity and bump functions. Related studies are still developing mainly by the author himself, and there exist related various explicit open problems.
\begin{Prob}
\label{prob:3}
For {\it infinite} graphs, formulate and answer to problems on realization of them as Reeb graphs of smooth functions of certain classes.
\end{Prob}
Here, an {\it infinite} graph means a $1$-dimensional locally finite CW complex which is not finite and the closure of each $1$-cell of which is homeomorphic to $D^1$.
In this note, we give an example as an answer to Problem \ref{prob:3}, respecting Problem \ref{prob:2} and our related studies (Theorem \ref{thm:1}). A map $c:X \rightarrow Y$ between topological spaces is said to be {\it proper} if the preimage $c^{-1}(K)$ of any compact set $K$ is compact.
\begin{Thm}
\label{thm:1}

There exists a connected infinite graph $G_1$ and for each integer $m>1$, we have an $m$-dimensional smooth connected manifold $X_{m} \subset {\mathbb{R}}^{m+1}$ which is non-compact and represented as the zero set ${e_m}^{-1}(0)$ of some real analytic function $e_m:{\mathbb{R}}^{m+1} \rightarrow \mathbb{R}$ and a function $c_{m}:X_{m} \rightarrow \mathbb{R}$ with the following properties.
\begin{enumerate}
\item \label{thm:1.1} The function $c_{m}$ is the restriction of ${\pi}_{m+1,1}$ to $X_{m}$. 
\item \label{thm:1.2} The Reeb space $R_{c_{m}}$ is also regarded as the Reeb graph of $c$ and two graphs $G_1$ and $R_{c_{m}}$ are isomorphic. 
\item \label{thm:1.3} The map $q_{c_{m}}:X_m \rightarrow R_{c_{m}}$ is proper and each regular contour of $c_{m}$ is diffeomorphic to $S^{m-1}$.
\end{enumerate}
\end{Thm}
\begin{Prob}
\label{prob:4}
Can we construct a compact manifold $X$ and a function $c:X \rightarrow \mathbb{R}$ with properties similar to ones in Problem \ref{prob:3}?
\end{Prob}
Since $X$ is compact, $R_c$ is compact, in Problem \ref{prob:4}. According to \cite[Theorem 7.5]{gelbukh1}, for a smooth function $c:X \rightarrow \mathbb{R}$ on a closed and connected manifold $X$, the Reeb space $R_c$ may not be a graph. However, $R_c$ is always homeomorphic to a so-called {\it Peano continuum} and homotopic to a (finite) graph: a {\it Peano continuum} is a separable, compact, connected and locally connected metrizable space. 
We also introduce the class of {\it graphs with ends} or {\it E-graphs}. If a CW complex is locally finite and the closure of each $1$-cell is homeomorphic to $D^1$ or $\{t \geq 0 \mid t \in \mathbb{R}\}$, then we call this an {\it E-graph}.
Theorem \ref{thm:2} is an answer to this and another new result. 
\begin{Thm}
\label{thm:2}

There exists a Peano continuum $G_2$ with the following properties.
\begin{enumerate}
\setcounter{enumi}{3}
\item \label{thm:2.1} The space $G_2$ is not homeomorphic to any E-graph and by choosing a point $p_{G_2} \in G_2$, being uniquely determined, $G_2-\{p_{G_2}\}$ is homeomorphic to and has the structure of an E-graph.
\item \label{thm:2.2} For any pair $(m_1,m_2)$ of integers greater than $1$, we have a smooth manifold $X_{m_1,m_2} \subset {\mathbb{R}}^{m+2}$ of dimension $m:=m_1+m_2$ being compact, connected and represented by $X_{m_1,m_2}={e_{m_1,m_2}}^{-1}(0)$ for some smooth map $e_{m_1,m_2}:{\mathbb{R}}^{m+2} \rightarrow {\mathbb{R}}^2$ being real analytic outside some set $Z_{m_1,m_2} \subset {\mathbb{R}}^{m+2}$ of Lebesgue measure $0$ and a function $c_{m_1,m_2}:X_{m_1,m_2} \rightarrow \mathbb{R}$ with the following properties.
\begin{enumerate}
\item \label{thm:2.2.1} The function $c_{m_1,m_2}$ is the restriction of ${\pi}_{m+2,1}$ to $X_{m_1,m_2}$. 
\item \label{thm:2.2.2} The Reeb space $R_{c_{m_1,m_2}}$ is homeomorphic to $G_2$. 
\item \label{thm:2.2.3} By choosing a point $p_{X_{m_1,m_2}} \in X_{m_1,m_2}$ suitably and considering the restriction of $c_{m_1,m_2}$ to $X_{m_1,m_2}-\{p_{X_{m_1,m_2}}\}$, we can induce the structure of an E-graph on $q_{c_{m_1,m_2}}(X_{m_1,m_2}-\{p_{X_{m_1,m_2}}\})$, where its vertex is defined by the rule applied for Reeb graphs of smooth functions. The resulting E-graph is also isomorphic to the E-graph $G_2-\{p_{G_2}\}$.

\end{enumerate}
\end{enumerate}
\end{Thm}
In the second section, we prove these theorems. Our main ingredient in the proof of Theorem \ref{thm:1} is a method first presented by the author in \cite{kitazawa2}, a method to construct real algebraic maps whose images are given regions surrounded by mutually disjoint real algebraic hypersurfaces in ${\mathbb{R}}^n$. This has been improved in the preprint \cite{kitazawa5} and is explained in a self-contained way and applied in proving Theorem \ref{thm:2}.
\section{A proof of Theorems \ref{thm:1} and \ref{thm:2}.}

\begin{proof}[A proof of Theorem \ref{thm:1}]

We first construct a manifold $X_m$ and a function $c_m:X_m \rightarrow \mathbb{R}$, in STEP 1-1. After that, we prove that the Reeb space $R_{c_m}$ is also the Reeb graph of $c_m$ and isomorphic to some $G_1$ for any $m>1$ and that the properties (\ref{thm:1.1}, \ref{thm:1.2}) are enjoyed, in STEP 1-2. In STEP 1-3, we prove that the property (\ref{thm:1.3}) is enjoyed to complete the proof. \\
\ \\
STEP 1-1 A manifold $X_m$ and a function $c_m:X_m \rightarrow \mathbb{R}$. \\
Let $c_{1,0}(x):=\frac{{\sin}^2  x}{2(x^2+1)}$ ($x \in \mathbb{R}$) and $c_{2,0}(x):=\frac{1}{x^2+1}$ ($x \in \mathbb{R}$).
We define $S_i:=\{(c_{i,0}(x),x) \mid x \in \mathbb{R}\}$ with $i=1,2$. We can see that the sets $S_1 \in {\mathbb{R}}^2$ and $S_2 \in {\mathbb{R}}^2$ are disjoint.
We consider the region $D_S:=\{(x_1,x_2) \mid x_2 \in \mathbb{R}, c_{1,0}(x_2)<x_1<c_{2,0}(x_2)\}$ and the closure $\overline{D_S}:=\{(x_1,x_2) \mid x_2 \in \mathbb{R}, c_{1,0}(x_2) \leq x_1 \leq c_{2,0}(x_2)\}$.

Respecting the original proof of \cite[Main Theorem 1]{kitazawa3}, we define
$X_{D_S,m}:=\{(x_1,x_2,y)=(x_1,x_2,{(y_j)}_{j=1}^{m-1}) \in {\mathbb{R}}^{m+1} \mid (x_1-c_{1,0}(x_2))(c_{2,0}(x_2)-x_1)-{\Sigma}_{j=1}^{m-1} {y_j}^2=0\} \subset {\mathbb{R}}^{m+1}$
and prove that this is a smooth manifold of dimension $m$ and the zero set of the real analytic function $e_m:{\mathbb{R}}^{m+1} \rightarrow \mathbb{R}$ defined canonically from $(x_1-c_{1,0}(x_2))(c_{2,0}(x_2)-x_1)-{\Sigma}_{j=1}^{m-1} {y_j}^2$. We also have the representation $X_{D_S,m}:=\{(x_1,x_2,y)=(x_1,x_2,{(y_j)}_{j=1}^{m-1}) \in \overline{D_S} \times {\mathbb{R}}^{m-1} \mid (x_1-c_{1,0}(x_2))(c_{2,0}(x_2)-x_1)-{\Sigma}_{j=1}^{m-1} {y_j}^2=0\}$. At a point $(x_{0,1},x_{0,2},y_0) \in X_{D_S}$ such that $(x_{0,1},x_{0,2}) \in D_S$, the partial derivative of $(x_1-c_{1,0}(x_2))(c_{2,0}(x_2)-x_1)-{\Sigma}_{j=1}^{m-1} {y_j}^2$ by some $y_j$ is not zero. At a point $(x_{0,1},x_{0,2},y_0) \in X_{D_S}$ such that $(x_{0,1},x_{0,2}) \in \overline{D_S}-D_S$, the partial derivative of $(x_1-c_{1,0}(x_2))(c_{2,0}(x_2)-x_1)-{\Sigma}_{j=1}^{m-1} {y_j}^2$ by $x_1$ is $c_{1,0}(x_{0,2})-x_{0,1}$ in the case $c_{2,0}(x_{0,2})-x_{0,1}=0$ ($c_{2,0}(x_{0,2})-x_{0,1}$ in the case $x_{0,1}-c_{1,0}(x_{0,2})=0$) and a non-zero real number. By the implicit function theorem, $X_{m}:=X_{D_S,m}$ is a real analytic manifold. We define a function $c_{m}:X_{m} \rightarrow \mathbb{R}$ as the restriction of ${\pi}_{m+1,1}$ to $X$ and prove that it is a desired function. \\
\ \\
STEP 1-2 The Reeb space $R_{c_m}$ is regarded as the Reeb graph of $c_m$ and isomorphic to some graph $G_1$, and the properties (\ref{thm:1.1}, \ref{thm:1.2}) are enjoyed. \\
We can easily see that the image of $c_{m}$ is $\{t \mid t \geq 0\}$ and that the restriction
of $c_{m}$ to ${c_{m}}^{-1}(\{t \mid t>0\})$ is proper. A point of $X_{m}$ with the form $(x_1,x_2,y) \in \mathbb{R} \times \mathbb{R} \times {\mathbb{R}}^{m-1}={\mathbb{R}}^{m+1}$ ($x_1 \neq 0$) is a critical point of $c$ if and only if it is of the form $(x_{c,1,1},x_{c,1,2},0) \in {\mathbb{R}}^2 \times {\mathbb{R}}^{m-1}$ with $x_{c,1,1}=c_{1,0}(x_{c,1,2})$ and $x_{c,1,2}$ being a critical point of the function $c_{m}$, or $(1,0,0) \in \mathbb{R} \times \mathbb{R} \times {\mathbb{R}}^{m-1}$. %, with $0 \in {\mathbb{R}}^{m-1}$ being the origin of course. 
Under the present situation, the points $(x_{c,1,1},x_{c,1,2})$ must satisfy the following.
\begin{itemize}
	\item $0<x_{c,1,2}<1$.
\item These points $x_{c,1,1}$ and $(x_{c,1,1},x_{c,1,2})$ appear discretely in $\mathbb{R}-\{0\}$ and ${\mathbb{R}}^2$, respectively, form closed subsets in these spaces, respectively, by the behavior of the function $c_0$. For a fixed non-zero number $x_{c,1,1}$ here, the set of all corresponding numbers $x_{c,1,2}$ must be finite.

\item Of course there exist infinitely many such points $(x_{c,1,1},x_{c,1,2}) \in {\mathbb{R}}^2$. 
\end{itemize}
\ \\

From this together with Ehresmann fibration theorem or more generally, \cite[Lemma 6.1]{gelbukh1}, stating that the Reeb space of a continuous function on compact and connected space whose level sets are connected is homeomorphic to $D^1$ or a one-point set, and Saeki's study on graph structures of Reeb spaces of smooth functions on compact manifolds such as \cite[Theorem 3.1]{saeki2} and \cite[Theorems 2.1 and 2.8 and "2 Reeb space and its graph structure"]{saeki3}, the Reeb space of the restriction of $c_{m}:X_{m} \rightarrow \mathbb{R}$ to ${c_{m}}^{-1}(\{t \mid t>0\})$ has the structure of an E-graph by defining vertices as in the definition of the Reeb graph of a smooth function. Let $c_{m,0}$ denote the restriction. More precisely, the following hold, where the restriction $c_{X_{m}}:X_{m} \rightarrow {\mathbb{R}}^2$ of ${\pi}_{m+1,2}$ to $X_{m}$ is used.
\begin{itemize}
\item Each contour of $c_{m}$ appears as the preimage ${c_{X_{m}}}^{-1}(I_P)$, where $I_P$ is a suitably chosen connected component of the intersection of $\overline{D_S}$ and a line $\{(p,t) \mid t \in \mathbb{R}\}$ with a suitably chosen value $p$ of $c_{m}$. There also exists a one-to-one correspondence between contours of $c_{m}$ and such connected components $I_P$.
\item Each regular (critical) contour of $c_{m}$ appears as the preimage ${c_{X_{m}}}^{-1}(I_{t_R})$ (resp. ${c_{X_{m}}}^{-1}(I_{t_C})$), where $I_{t_R}$ (resp. $I_{t_C}$) is a suitably chosen connected component of the intersection of $\overline{D_S}$ and a line $\{(t_R,t) \mid t \in \mathbb{R}\}$ (resp. $\{(t_C,t) \mid t \in \mathbb{R}\}$) with a suitably chosen regular (resp. critical) value $t_R$ (resp. $t_C$) of $c$. There also exists a one-to-one correspondence between regular (resp. critical) contours of $c_{m}$ and such connected components $I_{t_R}$ (resp. $I_{t_C}$). We call such a connected component $I_{t_R}$ (resp. $I_{t_C}$) a {\it regular} (resp. {\it critical}) {\it contour} of $\overline{D_S}$.
\item Each edge of $R_{c_{m,0}}$ appears as the image $q_{c_{m}}({c_{X_{m}}}^{-1}(I_{t_{C,1},t_{C,2}}))$, where $I_{t_{C,1},t_{C,2}}$ is defined as a suitably chosen connected component of $\overline{D_S} \bigcap \{(t_1,t_2) \mid t_{C,1}<t_1<t_{C,2}, t_2 \in \mathbb{R}\}$ for a suitably chosen pair $(t_{C,1},t_{C,2})$ of critical values of $c_{m}$ with the following constraints being posed.
\begin{itemize}
\item $0 \leq t_{C,1}<t_{C,2}\leq 1$.
\item The connected component $I_{t_{C,1},t_{C,2}}$ of $\overline{D_S} \bigcap \{(t_1,t_2) \mid t_{C,1}<t_1<t_{C,2}, t_2 \in \mathbb{R}\}$ is a disjoint union of regular contours of $\overline{D_S}$.
\item The number of critical contours of $\overline{D_S}$ whose intersections with the closure of $I_{t_{C,1},t_{C,2}}$ chosen in ${\mathbb{R}}^2$ are non-empty is $2$. One is of the form $I_{t_{C,1}}$ and the other is of the form $I_{t_{C,2}}$. Note that in the case $t_{C,1}=0$, ${c_{X_{m}}}^{-1}(I_{t_{C,1}})$ is not a contour of $c_{m,0}$, the restriction of $c_m$ to ${c_{m}}^{-1}(\{t \mid t>0\})$.
\end{itemize}
Furthermore, there also exists a one-to-one correspondence between these edges of the E-graph and such connected components $I_{t_{C,1},t_{C,2}}$. This naturally and uniquely defines the structure of an E-graph on $R_{c_{m,0}}$ and this does not depend on $m>1$, up to isomorphisms. 
\end{itemize}

Furthermore, ${c_{m}}^{-1}(0)$ consists of infinitely many points and is a discrete set in $X_{m}$. We can also see that the restriction of $q_{c_{m}}$ to ${c_{m}}^{-1}(0)$ is injective and by this ${c_{m}}^{-1}(0)$ is mapped onto ${\bar{c_{m}}}^{-1}(0) \subset R_{c_{m}}$. The set ${\bar{c_{m}}}^{-1}(0) \subset R_{c_{m}}$ is also regarded as a discrete set in $R_{c_{m}}$ and for each point $v_0$ of ${\bar{c_{m}}}^{-1}(0) \subset R_{c_{m}}$, we have its open neighborhood $N_{v_0}$ in $R_{c_{m}}$ homeomorphic to $\{t \mid 0 \leq t < 1\}$, containing no vertex of the E-graph explained just before, and making the restriction of $q_{c_m}$ to the preimage ${q_{c_m}}^{-1}(N_{v_0})$ proper. Thus $R_{c_{m}}$ is also seen as the Reeb graph of $c_{m}$. The property (\ref{thm:1.2}) has been shown to be enjoyed with the property (\ref{thm:1.1}). 

We have shown the fact that the Reeb graph $R_{c_{m}}$ is isomorphic to some infinite graph $G_1$ for any $m>1$.\\
\ \\
STEP 1-3 The property (\ref{thm:1.3}). \\
We prove that the property (\ref{thm:1.3}) is enjoyed.
By the fact that the restriction $c_{m,0}$ of $c_{m}$ to ${c_{m}}^{-1}(\{t \mid t>0\})$ is proper and the observation on ${c_{m}}^{-1}(0)$, ${\bar{c_{m}}}^{-1}(0)$, and a neighborhood $N_{v_0}$ of $v_0 \in{\bar{c_{m}}}^{-1}(0)$ in $R_{c_{m}}$, presented above, we can see that $q_{c_{m}}$ is also proper. Each regular contour of $c_{m}$ is always of the form ${c_{X_{m}}}^{-1}(I_R)$, where $I_R$ is a regular contour of $\overline{D_S}$ and diffeomorphic to $D^1$. From the structure of the maps and the manifolds, we can see that a regular contour of $c_{m}$ is diffeomorphic to $S^{m-1}$. \\
\ \\
This completes the proof of Theorem \ref{thm:1}.
\end{proof}
\begin{proof}[A proof of Theorem \ref{thm:2}]

We first construct a manifold $X_{m_1,m_2}$ and a function $c_{m_1,m_2}:X_{m_1,m_2} \rightarrow \mathbb{R}$ in STEP 2-1. After that, in STEP 2-2, we prove that the Reeb space $R_{c_{m_1,m_2}}$ is homeomorphic to some Peano continuum $G_2$ for any $m_1>1$ and $m_2>1$ to complete the proof. \\
\ \\
STEP 2-1 A manifold $X_{m_1,m_2}$ and a function $c_{m_1,m_2}:X_{m_1,m_2} \rightarrow \mathbb{R}$. \\
Let $c_0(x):=e^{-\frac{1}{x^2}} x {\sin}^2  (\frac{1}{x})$ ($x \in \mathbb{R}-\{0\}$).
For each positive integer $i$, the value of the $i$-th derivative ${c_0}^{(i)}$ of $c_0$ at $p_0$ is, by a simple calculation, shown to have the form $e^{-\frac{1}{{p_0}^2}} \times  {\Sigma}_{i_j \in I_i} \{p_{i_j}(\cos (\frac{1}{p_0}), \sin(\frac{1}{p_0}) ){(\frac{1}{p_0})}^{i_j}\} \times {(\frac{1}{p_0})}^{b_i}$, where the notation here is as follows.
\begin{itemize}
\item The set $I_i$ is a suitable non-empty finite set of non-negative integers. 
\item The function $p_{i_j}$ is a suitably chosen real polynomial function.
\item The number $b_i$ is a suitably chosen integer.
\end{itemize} 
The value $e^{-\frac{1}{{p_0}^2}} \times  {\Sigma}_{i_j \in I_i} \{p_{i_j}(\cos (\frac{1}{p_0}), \sin(\frac{1}{p_0}) ){(\frac{1}{p_0})}^{i_j}\} \times {(\frac{1}{p_0})}^{b_i}$ converges to $0$ as $p_0$ converges to $0$. The value $\frac{{c_0}^{(i)}(p_0)}{p_0}$ also converges to $0$ as $p_0$ converges to $0$. This implies that $c_0$ is extended to a continuous function on $\mathbb{R}$ uniquely. The resulting function $c_0:\mathbb{R} \rightarrow \mathbb{R}$ is real analytic outside $\{0\}$ and smooth. Note that $c_0(p_0)>0$ for $p_0>0$ that $c_0(p_0)<0$ for $p_0<0$, and that $c_0(0)=0$.

We define $S_1:=\{(c_0(x),x) \mid x \in \mathbb{R}\}$. We can also choose an ellipsoid $S_2:=\{(x_1,x_2) \in {\mathbb{R}}^2 \mid \frac{{(x_1-p_1)}^2}{r_1}+\frac{{(x_2-p_2)}^2}{r_2}=r_{1,2}\} \subset {\mathbb{R}}^2$ ($p_1>0$, $p_2>0$, $r_1>0$, $r_2>0$ and $r_{1,2}>0$) with the following properties. Let $D_{S_2}:=\{(x_1,x_2) \in {\mathbb{R}}^2 \mid \frac{{(x_1-p_1)}^2}{r_1}+\frac{{(x_2-p_2)}^2}{r_2}<r_{1,2}\} \subset {\mathbb{R}}^2$ and let $\overline{D_{S_2}} \subset {\mathbb{R}}^2$ denote the closure $\{(x_1,x_2) \in {\mathbb{R}}^2 \mid \frac{{(x_1-p_1)}^2}{r_1}+\frac{{(x_2-p_2)}^2}{r_2} \leq r_{1,2}\}$.
\begin{enumerate}
\setcounter{enumi}{5}
\item \label{thm:2.3} The set $S_1 \bigcap S_2$ is a two-point set. One point in this is $(0,0)$ and the other is represented as $(c_0(T),T)$ with a suitably chosen $T>0$. Furthermore, for the chosen $T$ and $\frac{{(x_1-p_1)}^2}{r_1}+\frac{{(x_2-p_2)}^2}{r_2}=r_{1,2}$, $T>p_2$ and $D_{S_2} \bigcap \{(c_0(t),t) \mid t \in \mathbb{R}\}=\{(c_0(t),t) \mid 0<t<T\}$. 

%The resulting ellipsoid is $S_2$.  
\item \label{thm:2.4} Furthermore, at each of these two points of $S_1 \bigcap S_2$ above, the tangent vector of the curve $S_1 \in {\mathbb{R}}^2$ and that of the curve $S_2 \in {\mathbb{R}}^2$ are always mutually independent if they are not zero vectors: a similar fact holds for their normal vectors at each of these two points above.

\item \label{thm:2.5} The set $D_{S_{1,2}}:=\{(x_1,x_2) \mid x_1 \in \mathbb{R}, c_0(x_2)<x_1\} \bigcap D_{S_2}$ and the closure $\overline{D_{S_{1,2}}}:=\{(x_1,x_2) \mid x_1, x_2 \in \mathbb{R}, c_0(x_2) \leq x_1\} \bigcap \overline{D_{S_2}}$ are non-empty. 
\end{enumerate}

We explain the situation enjoying these three properties, more precisely. We can choose $T:=T_0>0$ in such a way that $c_0(T_0) \geq c_0(t)>0$ for any $0 \leq t \leq T_0$ by considering ${\sin}^2 (\frac{1}{T_0})=1$. We first choose $S_{2,0}:=\{(x_1,x_2) \in {\mathbb{R}}^2 \mid \frac{{(x_1-c_0(T_0))}^2}{r_1}+{(x_2-\frac{T_0}{2})}^2=\frac{{T_0}^2}{4}\} \subset {\mathbb{R}}^2$ with $r_1>0$ being sufficiently large. After that, for a suitable number $s>1$, each point $(x_{1,0},x_{2,1}) \in S_{2,0}$ is moved to $(x_{1,0},x_{2,2})$ with the relation $T_0-x_{2,2}=s(T_0-x_{2,1})$, where the uniquely chosen point $(0,x_{2,0,1}) \in S_{2,0}$ of the two-point set $\{(0,x_{2,0,1}),(0,x_{2,0,2})\}=S_{2,0} \bigcap {{\pi}_{2,1}}^{-1}(0)$ ($0<x_{2,0,1}<x_{2,0,2}<T_0$) is moved to $(0,0)$ by this transformation. As $r_1$ goes to the infinity $\infty$, $s$ converges to $1$. By choosing a sufficiently large $r_1>0$ beforehand and a sufficiently small $s>1$, we have $S_2$ as the result of the transformation of $S_{2,0}$, with $p_1:=c_0(T_0)$. 

Respecting the original proof of \cite[Main Theorem 1]{kitazawa2}, we define
$X_{D_{S_{1,2}},m_1,m_2}:=\{(x_1,x_2,y_{S_1},y_{S_2})=(x_1,x_2,{(y_{S_1,j_1})}_{j_1=1}^{m_1},{(y_{S_1,j_2})}_{j_2=1}^{m_2}) \in {\mathbb{R}}^{m_1+m_2+2} \mid (x_1-c_0(x_2))-{\Sigma}_{j_1=1}^{m_1} {y_{S_1,j_1}}^2=0, r_{1,2}-\{\frac{{(x_1-p_1)}^2}{r_1}+\frac{{(x_2-p_2)}^2}{r_2}\}-{\Sigma}_{j_2=1}^{m_2} {y_{S_2,j_2}}^2=0\} \subset {\mathbb{R}}^{m_1+m_2+2}$ for integers $m_1>1$ and $m_2>1$, and we prove that this is a smooth, closed and connected manifold of dimension $m:=m_1+m_2$ and the zero set of the smooth map $e_{m_1,m_2}:{\mathbb{R}}^{m+2} \rightarrow {\mathbb{R}}^2$ defined canonically from $((x_1-c_0(x_2))-{\Sigma}_{j_1=1}^{m_1} {y_{S_1,j_1}}^2, r_{1,2}-\{\frac{{(x_1-p_1)}^2}{r_1}+\frac{{(x_2-p_2)}^2}{r_2}\}-{\Sigma}_{j_2=1}^{m_2} {y_{S_2,j_2}}^2)$. 
From (\ref{thm:2.5}), $X_{D_{S_{1,2}},m_1,m_2}$ is non-empty and we also have the representation\\ $X_{D_{S_{1,2}},m_1,m_2}:=\{(x_1,x_2,{(y_{S_1,j_1})}_{j_1=1}^{m_1}, {(y_{S_1,j_2})}_{j_2=1}^{m_2}) \in \overline{D_{S_{1,2}}} \times {\mathbb{R}}^{m} \mid (x_1-c_0(x_2))-{\Sigma}_{j_1=1}^{m_1} {y_{S_1,j_1}}^2=0, r_{1,2}-\{\frac{{(x_1-p_1)}^2}{r_1}+\frac{{(x_2-p_2)}^2}{r_2}\}-{\Sigma}_{j_2=1}^{m_2} {y_{S_2,j_2}}^2=0\}$. 

We investigate the rank of the differential of the map $e_{m_1,m_2}$. Our discussion here is done in a self-contained way, although readers may also refer to the preprint \cite{kitazawa5} of the author (\cite["Main Theorem 1" and "A proof of Main Theorem 1 with Theorem 1"]{kitazawa5}). \\
\ \\
Case 2-1-1 We consider a point $(x_{0,1},x_{0,2},y_{S_1,0},y_{S_2,0}) \in X_{D_{S_{1,2}}}$ such that $(x_{0,1},x_{0,2}) \in D_{S_{1,2}}$. There, the value of the partial derivative of $(x_1-c_0(x_2))-{\Sigma}_{j_1=1}^{m_1} {y_{S_1,j_1}}^2$ by some $y_{S_1,j_1}$ is not $0$ and that of the partial derivative of $r_{1,2}-\{\frac{{(x_1-p_1)}^2}{r_1}+\frac{{(x_2-p_2)}^2}{r_2}\}-{\Sigma}_{j_2=1}^{m_2} {y_{S_2,j_2}}^2$ by some $y_{S_2,j_2}$ is not $0$. The values of the partial derivative of $(x_1-c_0(x_2))-{\Sigma}_{j_1=1}^{m_1} {y_{S_1,j_1}}^2$ by each $y_{S_2,j_2}$ and the partial derivative of $r_{1,2}-\{\frac{{(x_1-p_1)}^2}{r_1}+\frac{{(x_2-p_2)}^2}{r_2}\}-{\Sigma}_{j_2=1}^{m_2} {y_{S_2,j_2}}^2$ by each $y_{S_1,j_1}$ there are both $0$.
 From this, the rank of differential of the map $e_{m_1,m_2}$ is $2$ at the point.   \\
\ \\
Case 2-1-2 Consider a point $(x_{0,1},x_{0,2},y_{S_1,0},y_{S_2,0}) \in X_{D_{S_{1,2}},m_1,m_2}$ such that $(x_{0,1},x_{0,2})$ is in $S_1$ and not in $S_2$. There, the value of the partial derivative of $r_{1,2}-\{\frac{{(x_1-p_1)}^2}{r_1}+\frac{{(x_2-p_2)}^2}{r_2}\}-{\Sigma}_{j_2=1}^{m_2} {y_{S_2,j_2}}^2$ by some $y_{S_2,j_2}$ is not $0$. The value of the partial derivative of $(x_1-c_0(x_2))-{\Sigma}_{j_1=1}^{m_1} {y_{S_1,j_1}}^2$ by $x_1$ there is $1$. The value of the partial derivative of $(x_1-c_0(x_2))-{\Sigma}_{j_1=1}^{m_1} {y_{S_1,j_1}}^2$ by each $y_{S_i,j_i}$ ($i=1,2$) there is $0$. From this, the rank of the differential of the map $e_{m_1,m_2}$ is $2$ at the point.   \\
\ \\
Case 2-1-3 Consider a point $(x_{0,1},x_{0,2},y_{S_1,0},y_{S_2,0}) \in X_{D_{S_{1,2}},m_1,m_2}$ such that $(x_{0,1},x_{0,2})$ is in $S_2$ and not in $S_1$. There, the value of the partial derivative of $r_{1,2}-\{\frac{{(x_1-p_1)}^2}{r_1}+\frac{{(x_2-p_2)}^2}{r_2}\}-{\Sigma}_{j_2=1}^{m_2} {y_{S_2,j_2}}^2$ by at least one of $x_1$ and $x_2$ is not $0$. The value of the partial derivative of $(x_1-c_0(x_2))-{\Sigma}_{j_1=1}^{m_1} {y_{S_1,j_1}}^2$ by some $y_{S_1,j_1}$ there is not $0$. The value of the partial derivative of $r_{1,2}-\{\frac{{(x_1-p_1)}^2}{r_1}+\frac{{(x_2-p_2)}^2}{r_2}\}-{\Sigma}_{j_2=1}^{m_2} {y_{S_2,j_2}}^2$ by each $y_{S_i,j_i}$ ($i=1,2$) there is $0$.
From this, the rank of the differential of the map $e_{m1,m_2}$ is $2$ at the point.\\
\ \\
Case 2-1-4  Consider a point $(x_{0,1},x_{0,2},y_{S_1,0},y_{S_2,0}) \in X_{D_{S_{1,2}},m_1,m_2}$ such that $(x_{0,1},x_{0,2})$ is in $S_1 \bigcap S_2$. There, the partial derivative of $r_{1,2}-\{\frac{{(x_1-p_1)}^2}{r_1}+\frac{{(x_2-p_2)}^2}{r_2}\}-{\Sigma}_{j_2=1}^{m_2} {y_{S_2,j_2}}^2$ by some $x_i$ is not $0$, and that of the partial derivative of $(x_1-c_0(x_2))-{\Sigma}_{j_1=1}^{m_1} {y_{S_1,j_1}}^2$ by some $x_i$ is not $0$. The normal vectors of the two curves $S_1 \in {\mathbb{R}}^2$ and $S_2 \in {\mathbb{R}}^2$ there are always mutually independent (, where they are non-zero vectors, of course), from (\ref{thm:2.3}, \ref{thm:2.4}). From this, the rank of the differential of the map $e_{m_1,m_2}$ is $2$ at the point.   \\
\ \\
From Cases 2-1-1. 2-1-2, 2-1-3 and 2-1-4, the rank of the differential of the map $e_{m_1,m_2}$ into ${\mathbb{R}}^2$, defined canonically from these two polynomials $(x_1-c_0(x_2))-{\Sigma}_{j_1=1}^{m_1} {y_{S_1,j_1}}^2$ and $r_{1,2}-\{\frac{{(x_1-p_1)}^2}{r_1}+\frac{{(x_2-p_2)}^2}{r_2}\}-{\Sigma}_{j_2=1}^{m_2} {y_{S_2,j_2}}^2$ is $2$, at each point of $X_{D_{S_{1,2}},m_1,m_2}$.
By implicit function theorem, $X_{D_{S_1,S_2},m_1,m_2}$ is a smooth manifold and the zero set of the map $e_{m_1,m_2}$. Note that in \cite[Main Theorem 1]{kitazawa5} (\cite[A proof of Main Theorem 1 with Theorem 1]{kitazawa5}), the real polynomial case is considered. Note also that the region $D_{S_{1,2}}$ and its dimension are considered in a certain general situation, in the preprint of the author. 

We put $X_{m_1,m_2}:=X_{D_{S_{1,2}},m_1,m_2}$. We define the set $Z_{m_1,m_2}:=\{(x_1,0,y_{S_1},y_{S_2}) \mid x_1 \in \mathbb{R}, {(y_{S_1,j_1})}_{j_1=1}^{m_1},{(y_{S_1,j_2})}_{j_2=1}^{m_2}) \in {\mathbb{R}}^{m_1} \times {\mathbb{R}}^{m_2}\} \subset {\mathbb{R}}^{m+2}$ as a subset of Lebesgue measure $0$ in ${\mathbb{R}}^{m+2}$. We define a function $c_{m_1,m_2}:X_{m_1,m_2} \rightarrow \mathbb{R}$ as the restriction of ${\pi}_{m+2,1}$ to $X_{m_1,m_2}$ and prove that this is a desired case. \\
\ \\
STEP 2-2 The Reeb space $R_{c_{m_1,m_2}}$ is homeomorphic to some Peano continuum $G_2$ which is not homeomorphic to any E-graph.\\
We can easily see that the image of $c_{m_1,m_2}$ is diffeomorphic to $D^1$. We can also see that the minimum is $0$ and that the maximum is $p_1+\sqrt{r_1r_{1,2}}$.
From (\ref{thm:2.3}, \ref{thm:2.4}) with the behavior of the function $c_0$ and the structures of $S_2$ and $X_{m_1,m_2}$ for example,
we can see that $X_{m_1,m_2}$ contains the origin $p_{X_{m_1,m_2}}:=0 \in {\mathbb{R}}^{m+2}$, that ${c_{m_1,m_2}}^{-1}(0)-\{p_{X_{m_1,m_2}}\}$ is a discrete set in $X_{m_1,m_2}$, that any open neighborhood of $p_{X_{m_1,m_2}}$ in $X_{m_1,m_2}$ must contain infinitely many points from ${c_{m_1,m_2}}^{-1}(0)-\{p_{X_{m_1,m_2}}\}$, and that the level set ${c_{m_1,m_2}}^{-1}(p_1+\sqrt{r_1r_{1,2}})$ is diffeomorphic to $S^{m_1-1}$.

A point of $X_{m_1,m_2}$ with the form $(x_1,x_2,y) \in \mathbb{R} \times \mathbb{R} \times {\mathbb{R}}^{m}={\mathbb{R}}^{m+2}$ ($x_1 \neq 0, p_1+\sqrt{r_1r_{1,2}}$) is a critical point of $c_{m_1,m_2}$ if and only if it is of the form $(x_{c,2,1},x_{c,2,2},y) \in {\mathbb{R}}^2 \times {\mathbb{R}}^{m}$ with $x_{c,2,1}=c_0(x_{c,2,2})$, $x_{c,2,2}$ being a critical point of the function $c_0$, and $y$ being an arbitrary point satisfying $(x_{c,2,1},x_{c,2,2},y) \in X$. Under the present situation, the points $(x_{c,2,1},x_{c,2,2})$ must satisfy the following.
\begin{itemize}
\item These points $x_{c,2,1}$ and $(x_{c,2,1},x_{c,2,2})$ appear discretely in $\mathbb{R}-\{0\}$ and ${\mathbb{R}}^2-\{(0,0)\}$ and form closed subsets in these spaces, respectively, by the behavior of the function $c_0$. For a fixed non-zero number $x_{c,2,1}$ here, the set of all corresponding numbers $x_{c,2,2}$ must be finite.
\item Of course there exist infinitely many such points $(x_{c,2,1},x_{c,2,2}) \in {\mathbb{R}}^2$. Furthermore, from the set of all of such points, we can find a sequence converging to $(0,0)$ in ${\mathbb{R}}^2$.
\end{itemize}

We consider $X_{m_1,m_2,0}:=X_{m_1,m_2}-\{p_{X_{m_1,m_2}}\}$. We also have ${q_{c_{m_1,m_2}}}^{-1}(q_{c_{m_1,m_2}}(p_{X_{m_1,m_2}}))=\{p_{X_{m_1,m_2}}\} \subset {\mathbb{R}}^{m+2}$.
We can also observe that ${\bar{c_{m_1,m_2}}}^{-1}(0)-\{q_{c_{m_1,m_2}}(p_{X_{m_1,m_2}})\}$ is a discrete subset of $R_{c_{m_1,m_2}}$ and that for each $v_0$ of the points of ${\bar{c_{m_1,m_2}}}^{-1}(0)-\{q_{c_{m_1,m_2}}(p_{X_{m_1,m_2}})\}$, we have a small connected open neighborhood $N_{v_0}$ in $R_{c_{m_1,m_2}}$ homeomorphic to $\{t \mid 0 \leq t < 1\}$. We can also see that any open neighborhood of $q_{c_{m_1,m_2}}(p_{X_{m_1,m_2}})$ chosen in $R_{c_{m_1,m_2}}$ contains infinitely many points from ${\bar{c_{m_1,m_2}}}^{-1}(0)-\{q_{c_{m_1,m_2}}(p_X)\}$ and this is the reason why $R_{c_{m_1,m_2}}$ is not homeomorphic to any $E$-graph.

From this together with arguments on the graph structure essentially same as ones in STEP 1-2, the Reeb space of the restriction of $c_{m_1,m_2}:X_{m_1,m_2} \rightarrow \mathbb{R}$ to $X_{m_1,m_2,0}$ is identified with $R_{c_{m_1,m_2}}-\{q_{c_{m_1,m_2}}(0)\}$ and has the structure of an E-graph by defining vertices as in the definition of the Reeb graph of a smooth function. The E-graph structure does not depend on $m_1>1$ and $m_2>1$ up to isomorphisms.
The original space $R_{c_{m_1,m_2}}$ is regarded as a one-point compactification of the E-graph $R_{c_{m_1,m_2}}-\{q_{c_{m_1,m_2}}(0)\}$ and its topology is determined uniquely from $R_{c_{m_1,m_2}}-\{q_{c_{m_1,m_2}}(0)\}$. It is homeomorphic to some Peano continuum $G_2$ which is not homeomorphic to any E-graph. We can choose the desired point $p_{G_2} \in G_2$ uniquely. \\
\ \\
By our proof, we can also see that the properties (\ref{thm:2.1}, \ref{thm:2.2}) are enjoyed. This completes the proof of Theorem \ref{thm:2}.
\end{proof}
We close this paper by presenting several problems. With a little effort, we may solve them, where the author has no idea.

\begin{Prob}
\label{prob:5}
This comes from Theorem \ref{thm:2}. 
Let $m>1$ and $k>0$ be integers. Let $c:X \rightarrow \mathbb{R}$ be a real analytic function with the following properties.
\begin{itemize}
\item The space $X$ is an $m$-dimensional smooth compact real manifold $X \subset {\mathbb{R}}^{m+k}$ represented as the zero set ${e_0}^{-1}(0)$ of some smooth map $e_0:{\mathbb{R}}^{m+k} \rightarrow {\mathbb{R}}^k$ which is real analytic outside a subset $Z \subset {\mathbb{R}}^{m+k}$ of Lebesgue measure $0$.
\item The function $c$ is represented as the restriction of ${\pi}_{m+k,1}$ to $X$.
\end{itemize}
Can we choose some finite or countably infinite set $F_c \subset R_c$ such that $R_c-F_c$ is homeomorphic to an E-graph? Or, as a stronger statement, is $R_c-F_c$ equipped with the structure of an E-graph defined by the rule in Reeb graphs, applied for the restriction of $c$ to $X-{q_c}^{-1}(F_c)$?
\end{Prob}
\begin{Prob}
\label{prob:6}
This comes from Theorem \ref{thm:1}. Let $m>1$ and $k>0$ be integers. Can we find a real analytic function $c:X \rightarrow \mathbb{R}$ with the following properties?
\begin{itemize}
\item The space $X$ is an $m$-dimensional non-compact real analytic manifold $X \subset {\mathbb{R}}^{m+k}$ represented as the zero set ${e_0}^{-1}(0)$ of some real analytic map $e_0:{\mathbb{R}}^{m+k} \rightarrow {\mathbb{R}}^k$.
\item The function $c$ is represented as the restriction of ${\pi}_{m+k,1}$ to $X$.
\item The image of $c$ is diffeomorphic to $D^1$.
\item The quotient map $q_c:X \rightarrow R_c$ is proper. The Reeb space $R_c$ is not equipped with the structure of the Reeb graph. As a stronger statement, $R_c$ is not homeomorphic to any graph.
\end{itemize}

\end{Prob}
\begin{Rem}
In important articles such as \cite{gelbukh1, saeki2, saeki3}, important and interesting examples on smooth functions and their Reeb spaces are presented. In our present study, the following properties are important.
\begin{itemize}
\item The map $q_c:X \rightarrow R_c$ is proper, where we abuse the notation we have used.
\item In Problems \ref{prob:5} and \ref{prob:6}, each component of the smooth map $e$ is represented by some single elementary function outside some subset of Lebesgue measure $0$ in ${\mathbb{R}}^{m+k}$. 
\end{itemize}
\end{Rem}

 \section{Conflict of interest and Data availability.}
  \noindent {\bf Conflict of interest.} \\
 The author is a researcher at Osaka Central Advanced Mathematical Institute (OCAMI researcher). This is supported by MEXT Promotion of Distinctive Joint Research Center Program JPMXP0723833165. He thanks this, where he is not employed there. \\
  %Some of works by other researchers and this version may overlap in some of the contents due to the nature that our problems are natural in theory of Morse functions and applications to differential topology and that related mathematical studies are very fundamental and classical in some senses, for example. However the present version of our paper is presented independent of these work. \\
  %Saga Souhatsu Mathematical Seminar (http://inasa.ms.saga-u.ac.jp/Japanese/saga-souhatsu.html), inviting the author as a speaker, is funded and supported by JST Fusion Oriented REsearch for disruptive Science and Technology JPMJFR202U: the author was a speaker on 2024/7/12 supported by this project.\\
  \ \\
  {\bf Data availability.} \\
  No data other than the present article is generated.

\end{document}